\overfullrule=0pt
\centerline {\bf Multiple periodic solutions of Lagrangian systems of relativistic oscillators}\par
\bigskip
\bigskip
\centerline {BIAGIO RICCERI}\par
\bigskip
\bigskip
{\bf Abstract.} - Let $B_L$ the open ball in ${\bf R}^n$ centered at $0$, of radius $L$, and let
$\phi$ be a homeomorphism from $B_L$ onto ${\bf R}^n$ such that $\phi(0)=0$ and $\phi=\nabla\Phi$,
where the function $\Phi:\overline {B_L}\to ]-\infty,0]$ is continuous and strictly convex in $\overline {B_L}$, and of class $C^1$ in $B_L$.
Moreover, let $F:[0,T]\times {\bf R}^n\to {\bf R}$ be a function which is measurable in $[0,T]$, of class $C^1$ in
${\bf R}^n$ and such that $\nabla_xF$ satisfies the $L^1$-Carath\'eodory conditions. Set
$$K=\{u\in \hbox {\rm Lip}([0,T],{\bf R}^n) : |u'(t)|\leq L\hskip 5pt for\hskip 3pt a.e.\hskip 3pt t\in [0,T] , u(0)=u(T)\}\ ,$$
and define the functional $I:K\to {\bf R}$ by
$$I(u)=\int_0^T(\Phi(u'(t))+F(t,u(t)))dt$$
for all $u\in K$. In [1], Brezis and Mawhin proved that any global minimum of $I$ in $K$ is a solution of the problem
$$\cases{(\phi(u'))'=\nabla_xF(t,u) & in $[0,T]$\cr & \cr
u(0)=u(T)\ , \hskip 3pt u'(0)=u'(T)\ .\cr}$$
In the present paper, we provide a set of conditions under which the functional $I$ has at least two global minima in $K$. This seems to be
the first result of this kind. The main tool of our proof is the well-posedness result obtained in [3].
\par
\bigskip
\bigskip
\bigskip
\bigskip
{\bf 1. - Introduction}\par
\bigskip
As the reader can notice, the title of the present paper is, intentionally, almost identical to the one of [1].\par
\smallskip
Actually, it is our aim to show how to obtain the multiplicity of periodic solutions for the systems mentioned in the title making
a joint use of the theory developed by Brezis and Mawhin in [1] with that we developed in [3].\par
\smallskip
To be more precise, we now fix some notations that we will keep throughout the paper and recall the main result of [1].\par
\smallskip
$L, T$ are two fixed positive numbers. For each $r>0$, we set $B_r=\{x\in {\bf R}^n :|x|<r\}$ ($|\cdot|$ being the Euclidean norm on ${\bf R}^n$)
and $\overline {B_r}$ is the closure of $B_r$.\par
\smallskip
We denote by ${\cal A}$ the family of all homeomorphisms $\phi$ from $B_L$ onto ${\bf R}^n$ such that $\phi(0)=0$ and $\phi=\nabla\Phi$,
where the function $\Phi:\overline {B_L}\to ]-\infty,0]$ is continuous and strictly convex in $\overline {B_L}$, and of class $C^1$ in $B_L$.
Notice that $0$ is the unique global minimum of $\Phi$ in $\overline {B_L}$.
\par
\smallskip
We denote by ${\cal B}$ the family of all functions $F:[0,T]\times {\bf R}^n\to {\bf R}$ which are measurable in $[0,T]$, of class $C^1$ in
${\bf R}^n$ and such that $\nabla_xF$ is measurable in $[0,T]$ and, 
for each $r>0$, one has $\sup_{x\in B_r}|\nabla_x F(\cdot,x)|\in L^1([0,T])$.\par
\smallskip
Given $\phi\in {\cal A}$ and $F\in {\cal B}$, we consider the problem
$$\cases{(\phi(u'))'=\nabla_xF(t,u) & in $[0,T]$\cr & \cr
u(0)=u(T)\ , \hskip 3pt u'(0)=u'(T)\ .\cr}\eqno{(P_{\phi,F})}$$
A solution of this problem is any function $u:[0,T]\to {\bf R}^n$ of class $C^1$, with $u'([0,T])\subset
B_L$, $u(0)=u(T)$, $u'(0)=u'(T)$,  such that the composite function $\phi\circ u'$ is absolutely
continuous in $[0,T]$ and one has $(\phi\circ u')'(t)=\nabla_xF(t,u(t))$ for a.e. $t\in [0,T]$.\par
\smallskip
Now, we set
$$K=\{u\in \hbox {\rm Lip}([0,T],{\bf R}^n) : |u'(t)|\leq L\hskip 5pt for\hskip 3pt a.e.\hskip 3pt t\in [0,T] , u(0)=u(T)\}\ ,$$
Lip$([0,T],{\bf R}^n)$ being the space of all Lipschitzian functions from $[0,T]$ into ${\bf R}^n$.\par
\smallskip
Clearly, one has
$$\sup_{[0,T]}|u|\leq LT + \inf_{[0,T]}|u| \eqno{(1.1)}$$
for all $u\in K$. To see this, take $t_0\in [0,T]$ such that $|u(t_0)|=\inf_{[0,T]}|u|$ and observe that, for each $t\in [0,T]$, one has
$$|u(t)-u(t_0)|=\left | \int_{t_0}^tu'(\tau)d\tau\right |\leq LT\ .$$\par
\smallskip
Next, consider the functional $I:K\to {\bf R}$ defined by
$$I(u)=\int_0^T(\Phi(u'(t))+F(t,u(t)))dt$$
for all $u\in K$.\par
\smallskip
The basic result of the theory developed in [1] is as follows:\par
\medskip
THEOREM 1.1 ([1], Theorem 5.2). - {\it Any global minimum of $I$ in $K$ is a solution of problem $(P_{\phi,F})$.}\par
\medskip
Well, the aim of the present paper is to provide a set of conditions under which the functional $I$ has at least two global minima in $K$.\par
\smallskip
As far as we know, this is the first result of this kind, and so we cannot do any proper comparison with previous ones.\par
\smallskip
Notice that some multiplicity results for problem $(P_{\phi,F})$ are already available in the literature. In this connection, we refer to the numerous
references contained in the very recent survey by Mawhin [2] and in [4]. But, as we repeat, in those papers the multiple solutions of problem 
$(P_{\phi,F})$ are not shown to be global minima of the functional $I$ in $K$.\par
\smallskip
As we said at the beginning, our main tool is provided by the main result obtained in [3] which is recalled in the next section.\par
\smallskip
Finally, Section 3 contains the statement of our multiplicity result, its proof and various related remarks.\par
\bigskip
{\bf 2. - A well-posedness theorem}\par
\bigskip
In this section, we summarize the theory developed in [3].\par
\smallskip
So, let $X$ be a Hausdorff topological space, $J, \Psi$  two
real-valued functions defined in $X$, and $a, b$  two numbers in $[-\infty,+\infty]$,
with $a<b$. \par
\smallskip
If $a\in {\bf R}$ (resp. $b\in {\bf R}$), we denote by $M_a$ (resp. $M_b$) 
the set of all global minima of the function $J+a\Psi$ (resp. $J+b\Psi$), while
if $a=-\infty$ (resp. $b=+\infty$), $M_a$ (resp. $M_b$) stands for the empty set.
We adopt the conventions $\inf\emptyset=+\infty$, $\sup\emptyset=-\infty$.
\smallskip
We also set
$$\alpha=\max\left \{ \inf_X \Psi,\sup_{M_b}\Psi\right \}\ ,$$
$$\beta=\min\left \{ \sup_X \Psi,\inf_{M_a}\Psi\right \}\ .$$
One proves that $\alpha\leq \beta$.\par
\smallskip
As usual, given a function $f:X\to {\bf R}$ and a set $C\subseteq X$, 
  we say that
the problem of minimizing $f$ over $C$ is well-posed if the following two
conditions hold:\par
\smallskip
\noindent
-\hskip 7pt the restriction of $f$ to $C$ has a unique global minimum,
say $\hat x$\ ;\par
\smallskip
\noindent
-\hskip 7pt every sequence $\{x_n\}$ in $C$ such that $\lim_{n\to \infty}f(x_n)=
\inf_Cf$, converges to $\hat x$.\par
\smallskip
A set of the type $\{x\in X: f(x)\leq r\}$ is said to be a sub-level set of
$f$. 
\smallskip
The main result of [3] is as follows:\par
\medskip
THEOREM 2.1 ([3], Theorem 1). - {\it 
Assume that $\alpha<\beta$ and that, for each $\lambda\in
]a,b[$, the function $J+\lambda\Psi$ has sequentially
compact sub-level sets and admits a unique global minimum in
$X$.\par
Then, for each $r\in ]\alpha,\beta[$, 
the problem of minimizing $J$ over $\Psi^{-1}(r)$ is well-posed.\par
Moreover,
if we denote by $\hat x_r$  the unique global minimum
of $J_{|\Psi^{-1}(r)}$ $(r\in ]\alpha,\beta[)$,
the functions $r\to \hat x_r$ and $r\to J(\hat x_r)$ are continuous in
$]\alpha,\beta[$, and, for some $\hat \lambda_r\in ]a,b[$, $\hat x_r$ is the global minimum in $X$
of the function $J+\hat\lambda_r\Psi$ }\par 
\bigskip
{\bf 3. - The main result}\par
\bigskip
Here is our main result:\par
\medskip
THEOREM 3.1. - {\it Let $\phi\in {\cal A}$, $F\in {\cal B}$, $G\in C^{1}({\bf R}^n)$, $\psi\in L^{1}([0,T])\setminus \{0\}$, with $\psi\geq 0$.
 Moreover, let $\gamma:[0,+\infty[\to {\bf R}$ be a 
convex strictly increasing function such that $\lim_{s\to +\infty}{{\gamma(s)}\over {s}}=+\infty$.
Assume that the following assumptions are
satisfied:\par
\noindent
$(i_1)$\hskip 3pt 
for a.e. $t\in [0,T]$ and for every $x\in {\bf R}^n$, one has
$$\gamma(|x|)\leq F(t,x)\ ;$$
\noindent
$(i_2)$\hskip 3pt $\liminf_{|x|\to +\infty}{{G(x)}\over {|x|}}>-\infty\ ;$\par
\noindent
$(i_3)$\hskip 3pt the function $G$ has no global minima in ${\bf R}^n$\ ;\par
\noindent
$(i_4)$\hskip 3pt there exist two points $x_1, x_2\in {\bf R}^n$ such that
$$\inf_{x\in {\bf R}^n}\int_0^TF(t,x)dt<\max\left \{ \int_0^TF(t,x_1)dt, \int_0^TF(t,x_2)dt\right\}$$
and
$$G(x_1)=G(x_2)=\inf_{B_c}G$$
where
$$c=LT+\gamma^{-1}\left ( {{1}\over {T}}\max\left \{ \int_0^TF(t,x_1)dt, \int_0^TF(t,x_2)dt\right\}\right )\ .$$
Then, there exist $\tilde\lambda>0$ such that the problem
$$\cases{(\phi(u'))'=\nabla_x(F(t,u)+\tilde\lambda \psi(t)G(u)) & in $[0,T]$\cr & \cr
u(0)=u(T)\ ,\hskip 3pt u'(0)=u'(T)\cr}$$
has at least two solutions which are global minima in $K$ of the functional 
$$u\to \int_0^T(\Phi(u'(t))+F(t,u(t))+\tilde\lambda \psi(t)G(u(t)))dt\ .$$}\par
\smallskip
PROOF. Let $C^0([0,T], {\bf R}^n)$ be the space of all continuous functions from $[0,T]$ into ${\bf R}^n$, with the norm $\sup_{[0,T]}|u|$.
We are going to apply Theorem 2.1 taking: $a=0$, $b=+\infty$, $X=K$ regarded as a subset of $C^0([0,T], {\bf R}^n)$ with the relative topology
 and
$$J(u)=\int_0^T\psi(t)G(u(t))dt\ ,$$
$$\Psi(u)=\int_0^T(\Phi(u'(t))+F(t,u(t)))dt$$
for all $u\in K$. Fix $\lambda>0$. 
By $(i_2)$, for a suitable constant $\delta>0$, we have
$$-\delta(|x|+1)\leq G(x)$$
for all $x\in {\bf R}^n$. 
For each $u\in K$, in view of $(i_1)$, $(1.1)$ and of the convexity of $\gamma$,
using Jensen inequality, we get
$$\int_0^T\psi(t)G(u(t))dt+\lambda\int_0^T(\Phi(u'(t))+F(t,u(t)))dt\geq -\delta\int_0^T\psi(t)|u(t)|dt+\lambda\int_0^T\gamma(|u(t)|)dt
-\delta\int_0^T\psi(t)dt+\lambda\Phi(0)T\geq$$
$$-\delta\int_0^T\psi(t)dt\sup_{[0,T]}|u|+\lambda T\gamma\left ( {{1}\over {T}}\int_0^T|u(t)|dt\right )-\delta\int_0^T\psi(t)dt+\lambda\Phi(0)T\geq$$
$$-\delta\int_0^T\psi(t)dt\sup_{[0,T]}|u|+\lambda T\gamma\left (\inf_{[0,T]}|u|\right )-\delta\int_0^T\psi(t)dt+\lambda\Phi(0)T\geq$$
$$-\delta\int_0^T\psi(t)dt\sup_{[0,T]}|u|+\lambda T \gamma\left (\sup_{[0,T]}|u|-LT\right )-\delta\int_0^T\psi(t)dt+\lambda\Phi(0)T\ .\eqno{(3.1)}$$
In turn, since $\lim_{s\to +\infty}{{\gamma(s-LT)}\over {s}}=+\infty$, we infer from $(3.1)$ that, for every
$\rho\in {\bf R}$, there is $M>0$ such that 
$$\left \{u\in K : \int_0^T\psi(t)G(u(t))dt+\lambda\int_0^T(\Phi(u'(t))+F(t,u(t)))dt\leq\rho\right\}\subseteq
\left\{u\in K : \sup_{[0,T]}|u|\leq M\right\}\ .\eqno{(3.2)}$$
Now, observe that $K$ is a closed subset of $C^0([0,T], {\bf R}^n)$. On the other hand, from Lemma 4.1 of [1], it
follows that the functional $J+\lambda\Psi$ is lower semicontinuous in $K$. Summarizing: the functions belonging to any sub-level set of $J+\lambda\Psi$ are
equi-continuous (since they are in $K$) and equi-bounded in view of $(3.2)$. Hence, by the Ascoli-Arzel\`a theorem, any sub-level set of $J+\lambda\Psi$ is
relatively sequentially compact in $C^0([0,T],{\bf R}^n)$. But, for the remarks above, the same set is closed in $C^0([0,T],{\bf R}^n)$, and so it is
sequentially compact in $K$. Next, observe that the functional $J$ ha no global minima in $K$. Since the constant functions lie in $K$, it is clear that
$$\inf_KJ=\inf_{{\bf R}^n}G\int_0^T\psi(t)dt\ .$$
Hence, if $G$ is unbounded below, so $J$ is too. Now, suppose that $G$ is bounded below. Arguing by contradiction, assume that $\hat u\in K$ is a global
minimum of $J$. Then, we would have
$$\int_0^T\psi(t)\left (G(\hat u(t))-\inf_{{\bf R}^n}G\right )dt=0\ ,$$
and so, since the integrand is non-negative, it would follow
$$\psi(t)\left (G(\hat u(t))-\inf_{{\bf R}^n}G\right )=0$$
for a.e. $t\in [0,T]$. Therefore, since $\psi\neq 0$, for some $t\in [0,T]$, we would have $G(\hat u(t))=\inf_{{\bf R}^n}G$, against $(i_3)$.
Notice that the absence of global minima for $J$ implies that 
$$\beta=\sup_K\Psi\ .$$
Moreover, since $\lim_{s\to +\infty}\gamma(s)=+\infty$, from $(i_1)$ it follows that 
$$\sup_K\Psi=+\infty\ .$$ 
Furthermore, since $b=+\infty$, we have
$$\alpha=\inf_K\Psi\ .$$
Clearly
$$\inf_K\Psi\leq \inf_{x\in {\bf R}^n}\int_0^TF(t,x)dt+\Phi(0)T\ .$$
Now, put
$$r=\max\left \{ \int_0^TF(t,x_1)dt, \int_0^TF(t,x_2)dt\right\}+\Phi(0)T\ .$$
By the above remarks and by the inequality in $(i_4)$, we have
$$\alpha<r<\beta\ .$$
Fix $u\in \Psi^{-1}(]-\infty,r])$. By $(i_1)$ and Jensen inequality again, we have
$$r\geq\int_0^T(\Phi(u'(t))+F(t,u(t)))dt\geq \int_0^T\gamma(|u(t)|)dt+\Phi(0)T\geq T\gamma\left ({{1}\over {T}}\int_0^T|u(t)|dt\right )+
\Phi(0)T\ ,$$
and so
$$\gamma\left ({{1}\over {T}}\int_0^T|u(t)|dt\right )\leq {{1}\over {T}}\max\left \{ \int_0^TF(t,x_1)dt, \int_0^TF(t,x_2)dt\right\}\ .$$
Applying $\gamma^{-1}$, we get
$${{1}\over {T}}\int_0^T|u(t)|dt\leq \gamma^{-1}\left ( {{1}\over {T}}\max\left \{ \int_0^TF(t,x_1)dt, \int_0^TF(t,x_2)dt\right\}\right )$$
and hence
$$\inf_{[0,T]}|u|\leq \gamma^{-1}\left ( {{1}\over {T}}\max\left \{ \int_0^TF(t,x_1)dt, \int_0^TF(t,x_2)dt\right\}\right )\ .$$
In view of $(1.1)$, we then infer that
$$\sup_{[0,T]}|u|\leq LT+\gamma^{-1}\left ( {{1}\over {T}}\max\left \{ \int_0^TF(t,x_1)dt, \int_0^TF(t,x_2)dt\right\}\right )\ .\eqno{(3.3)}$$
In turn, in view of $(i_4)$, $(3.3)$ implies that
$$J(x_1)=J(x_2)\leq J(u)\ .$$
Since $x_1, x_2\in \Psi^{-1}(]-\infty,r])$, we then conclude that $x_1, x_2$ are two distinct global minima
of $J_{|\Psi^{-1}(]-\infty,r])}$. 
Now, arguing by contradiction, assume that, for every $\lambda>0$, the functional $J+\lambda\Psi$ has a unique
global minimum in $K$. 
Then, by Theorem 2.1 (recall that $J+\lambda\Psi$ has sequentially compact sub-level sets), there would exist
$\hat\lambda_r>0$ and $\hat u_r\in \Psi^{-1}(r)$ such that $\hat u_r$ is the unique global minimum of $J+\hat\lambda_r\Psi$ in $K$. Then, for
$i=1,2$, we would have
$$\inf_{u\in K}(J(u)+\hat\lambda_r\Psi(u))\leq J(x_i)+\hat\lambda_r\Psi(x_i)\leq J(\hat u_r)+\hat\lambda_r\Psi(\hat u_r)=$$
$$\inf_{u\in K}(J(u)+\hat\lambda_r\Psi(u))\ .$$
That is to say, $x_1$ and $x_2$ would be two distinct global minima in $K$ of the functional $J+\hat\lambda_r\Psi$, a contradiction. So, 
 there exists
some $\hat \lambda>0$ such that the functional $J+\hat\lambda\Psi$ has at least two global minima in $K$. To conclude the proof, take
$\tilde \lambda={{1}\over {\hat\lambda}}$ and apply Theorem 1.1.\hfill $\bigtriangleup$\par
\medskip
In the sequel, the following further result from [1] will be useful:\par
\medskip
PROPOSITION 3.1 ([1], Proposition 3.2). - {\it Let $\phi\in {\cal A}$, $p>1$ and $\mu>0$.\par
Then, for every $\omega\in L^1([0,T],{\bf R}^n)$, the problem
$$\cases{(\phi(u'))'=\mu |u|^{p-2}u+\omega(t) & in $[0,T]$\cr & \cr
u(0)=u(T)\ ,\hskip 3pt u'(0)=u'(T)\cr}$$
has a unique solution.}
\par
\medskip
The next three examples (where $\phi\in {\cal A}$) show that, in Theorem 3.1, none of $(i_2)-(i_4)$ can be removed at all.\par
\medskip
EXAMPLE 3.1. - Take: $F(x)={{|x|^2}\over {2}}$, $G(x)=\langle z,x\rangle$, with $z\in {\bf R}^n\setminus \{0\}$,
$\psi=1$, $\gamma(s)={{s^2}\over {2}}$. Clearly, $(i_1)-(i_3)$ are satisfied, but, for
every $\lambda\in {\bf R}$,  the problem
$$\cases{(\phi(u'))'=u+\lambda z & in $[0,T]$\cr & \cr
u(0)=u(T)\ ,\hskip 3pt u'(0)=u'(T)\cr}$$
has a unique solution by Proposition 3.1.\par
\medskip
EXAMPLE 3.2. - Take: $F(x)={{|x|^2}\over {2}}$, $G=0$, $\gamma(s)={{s^2}\over {2}}$. Clearly, $(i_1)$, $(i_2)$ and $(i_4)$ are
satisfied, but, by Proposition 3.1, $0$ is the unique solution of the problem
$$\cases{(\phi(u'))'=u & in $[0,T]$\cr & \cr
u(0)=u(T)\ ,\hskip 3pt u'(0)=u'(T)\ .\cr}$$
\medskip
EXAMPLE 3.3. - Take: $F(x)=|x|^2$, 
$$G(x)=\cases {0 & if $|x|\leq LT+1$\cr & \cr
-(|x|-LT-1)^3 & if $|x|>LT+1$\ ,\cr}$$
$\psi=1$, $\gamma(s)=s^2$. Clearly, $(i_1)$, $(i_3)$ and $(i_4)$ are satisfied. In particular, $(i_4)$ is satisfied by any pair of distinct
points $x_1, x_2\in {\bf R}^n$ such that $|x_1|=|x_2|\leq 1$. However, for each $\lambda>0$, the functional $u\to \int_0^T(\Phi(u'(t))+|u(t)|^2+
\lambda G(u(t)))dt$ is unbounded below in $K$.\par
\medskip
We conclude with a joint consequence of Theorem 3.1 and Proposition 3.1.\par
\medskip
THEOREM 3.2. - {\it Let $\phi\in {\cal A}$, $p>1$, $G\in C^1({\bf R}^n)$,
$\psi\in L^{1}([0,T])\setminus \{0\}$, with $\psi\geq 0$. Assume that $G$ satisfies assumptions $(i_2)$, $(i_3)$ of Theorem 3.1
and the following:\par
\noindent
$(j_1)$\hskip 3pt there is $\rho>LT$ such that $G$ is constant in $B_{\rho}$\ .\par
Then, there exists $\tilde\lambda>0$ such that the problem
$$\cases{(\phi(u'))'=|u|^{p-2}u+\tilde\lambda \psi(t)\nabla G(u) & in $[0,T]$\cr & \cr
u(0)=u(T)\ ,\hskip 3pt u'(0)=u'(T)\cr}$$
has at least one solution which is a global minimum in $K$ of the functional 
$$u\to \int_0^T\left (\Phi(u'(t))+{{|u(t)|^p}\over {p}}+\tilde\lambda \psi(t)G(u(t))\right )dt$$
and whose range is contained in ${\bf R}^n\setminus \overline {B_{\rho-LT}}$\ .}\par
\smallskip
PROOF. Apply Theorem 3.1 with $F(t,x)={{|x|^p}\over {p}}$ and $\gamma(s)={{s^p}\over {p}}$.
Concerning $(i_4)$, notice that it is satisfied by any pair of distinct
distinct points $x_1, x_2\in {\bf R}^n$, such
that $|x_1|=|x_2|\leq \rho-LT$. This comes from $(j_1)$ after observing that
$\gamma^{-1}\left ({{1}\over {T}}\int_0^TF(t,x_1)dt\right )=|x_1|$. So, by Theorem 3.1, there exists $\tilde\lambda>0$
such that  the problem
$$\cases{(\phi(u'))'=|u|^{p-2}u+\tilde\lambda \psi(t)\nabla G(u) & in $[0,T]$\cr & \cr
u(0)=u(T)\ ,\hskip 3pt u'(0)=u'(T)\cr}$$
has at least one non-zero solution which is a global minimum in $K$ of the functional 
$$u\to \int_0^T\left (\Phi(u'(t))+{{|u(t)|^p}\over {p}}+\tilde\lambda \psi(t)G(u(t))\right )dt\ .$$
Denote by $w$ such a solution. To complete the proof, we have to show that $\inf_{[0,T]}|w|>\rho-LT$.
Arguing by contradiction, assume that $\inf_{[0,T]}|w|\le\rho-LT$. Then, by $(1.1)$, we would have
$$\sup_{[0,T]}|w|\leq \rho\ .$$
By $(j_1)$, this would imply that $w$ is a solution of the problem
$$\cases{(\phi(u'))'=|u|^{p-2}u & in $[0,T]$\cr & \cr
u(0)=u(T)\ ,\hskip 3pt u'(0)=u'(T)\cr}$$
and hence $w=0$ by Proposition 3.1, which is a contradiction. The proof is complete.\hfill $\bigtriangleup$
\vfill\eject
\centerline {\bf References}\par
\bigskip
\bigskip
\noindent
[1]\hskip 5pt H. BREZIS and J. MAWHIN, {\it Periodic solutions of Lagrangian systems of relativistic oscillators},
Commun. Appl. Anal., {\bf 15} (2011), 235-250.\par
\medskip
\noindent
[2]\hskip 5pt J. MAWHIN, {\it Multiplicity of solutions of relativistic-type systems with periodic nonlinearities: a survey},
Electron. J. Differ. Equ. Conf., {\bf 23} (2016), 77-86.\par
\medskip
\noindent
[3]\hskip 5pt B. RICCERI, {\it Well-posedness of constrained minimization problems via saddle-points}, J. Global Optim., {\bf 40} (2008),
389-397.\par
\medskip
\noindent
[4]\hskip 5pt X. WANG, Q. LIU and D. QIAN, {\it Existence and multiplicity results for some nonlinear problems with singular 
$\phi$-Laplacian via a geometric approach},  Bound. Value Probl. {\bf 2016}, 2016:47.
\bigskip
\bigskip
\bigskip
\bigskip
Department of Mathematics\par
University of Catania\par
Viale A. Doria 6\par
95125 Catania, Italy\par
{\it e-mail address}: ricceri@dmi.unict.it

\bye